\documentclass[12pt]{amsart}
\usepackage{amsthm,amsopn,amsfonts,amssymb,pb-diagram}
\usepackage{hyperref}

\usepackage{euler, amsfonts, amssymb, latexsym, epsfig,epic}

\font\co=lcircle10
\def\boxcross{\ \smash{\lower6.5pt\hbox{\rlap{\hskip4.5pt\vrule height13.5pt}}
                \raise0pt\hbox{\rlap{\hskip-2pt \vrule height.4pt depth0pt
                        width13.5pt}}}\hskip12.7pt}

\def\boxelbow{\ \hskip.1pt\smash{%
               \hbox{\co \hskip 5.5pt\rlap{\mathsurround=0pt\rlap{\mathsurround=0pt\char'006}\lower0.4pt\rlap{\char'004}}
                \lower6.5pt\rlap{\hskip-0.2pt\vrule height3pt}
                \raise3.5pt\rlap{\hskip-0.2pt\vrule height3.2pt}}
                \hbox{%
                  \rlap{\hskip-6.4pt \vrule height.4pt depth0pt
width2.5pt}%
                  \rlap{\hskip4.05pt \vrule height.4pt depth0pt
width3.1pt}}}
                \hskip8.7pt}

\newtheorem{Theorem}{Theorem}[section]
\newtheorem*{Theorem*}{Theorem}

\newtheorem{Corollary}{Corollary}

\theoremstyle{remark}

\newcommand\defn[1]{{\bf #1}}

\begin{document}

\title{Compatibly Frobenius split subschemes are rigid}
\author{Allen Knutson}
\thanks{AK was supported by an NSF grant.}
\address{Department of Mathematics, Cornell University,
        Ithaca, NY 14853, USA}
\email{allenk@math.cornell.edu}
\date{January 15, 2009}

\maketitle

\begin{abstract}
  Schwede proved very recently that in a quasiprojective scheme $X$ with a
  fixed Frobenius splitting, there are only finitely many subschemes $\{Y\}$
  that are compatibly split. (A simpler proof has already since been given,
  by Kumar and Mehta.) It follows that their deformations 
  (as compatibly split subschemes) are obstructed.
  
  We give a short proof that if $X$ is projective, 
  its compatibly split subschemes $\{Y\}$ have no deformations at all
  (again, as compatibly split subschemes).
  This reproves Schwede's result in some simple cases.
\end{abstract}

\setcounter{tocdepth}{1}

\section{The rigidity theorem}

\newcommand\FF{{\mathbb F}}
\newcommand\naturals{{\mathbb N}}
Let $\FF$ be a perfect field of characteristic $p>0$, 
and $R$ a commutative Noetherian ring containing $\FF$.
A \defn{Frobenius splitting} of $R$ is a map $\varphi:R \to R$ such that
$$ \varphi(a+b) = \varphi(a) + \varphi(b),\quad
\varphi(a^p b) = a \varphi(b), \quad
\varphi(1) = 1. $$
Our reference for Frobenius splittings is \cite{BrionKumar}. 

An ideal $I\leq R$ is called \defn{compatibly split} if $\varphi(I)
\subseteq I$. In this case $\varphi(I) = I$ and $\varphi$ descends to
give a splitting of $R/I$. Though the definition has been around
for twenty years, only recently was the following fundamental
observation made:

\begin{Theorem*}\cite[Corollary to Theorem 4.5]{Schwede}\cite{KumarMehta}
  Let $R$ be a ring with a Frobenius splitting $\varphi$.
  Then there are only finitely many ideals $I$ that are compatibly split.
\end{Theorem*}

In a local ring, this follows from \cite[Corollary 3.2]{EH},
which says that {\em any} set of radical ideals closed under sum,
intersection, and primary decomposition is finite. These properties
are easy to show of the set of compatibly split ideals 
\cite[Proposition 1.2.1]{BrionKumar}. 
Rather than the local situation (also 
in \cite{Sharp}), we will study the 
related
graded situation.

In what follows $R = \bigoplus_{n\in\naturals} R_n$ will be a graded ring 
with $R_0 = \FF$. The \defn{Hilbert polynomial} of a graded ideal $I\leq R$
is the unique polynomial $h_I$ such that $h_I(n) = \dim_\FF (R/I)_n$
for all sufficiently divisible $n$. The \defn{Hilbert scheme} $H_f(R)$,
for which our reference is \cite{WhySchemes}, represents the 
functor of saturated graded ideals $I\leq R$ with $h_I = f$.

Call $\varphi$ a \defn{graded splitting} if $\varphi(R_k) = 0$ unless
$p|k$, in which case $\varphi(R_k) \subseteq R_{k/p}$.  Asking
$\varphi$ to be graded is not much of a hardship, as from any
splitting $\varphi$ one may extract the \defn{graded part} $\varphi'$,
by defining $\varphi'(r)$ for $r \in R_k$ as the $k/p$-homogeneous
component of $\varphi(r)$ if $p|k$, and $\varphi'(r)=0$ otherwise. 
(I thank Michel Brion for this observation.) The same trick works for
multigradings, where $k$ lives in a lattice.

\begin{Corollary}
  Let $R$ be a graded ring with a graded splitting and $R_0 = \FF$.
  Let $f$ be a Hilbert polynomial. Then the set of compatibly split
  saturated graded ideals forms a finite subset of $H_f(R)$.
\end{Corollary}

We will see this as a corollary of another theorem, Theorem 
\ref{thm:main} below. To do this we will define a natural sub{\em scheme}
$H_f(R,\varphi)$ of $H_f(R)$ supported on the set of compatibly split ideals,
and show that its Zariski tangent spaces vanish. 

First we recall Mumford's closed embedding of the Hilbert scheme:

\begin{Theorem*}
  \label{Mumford}
  Fix $R,f$. Then $\exists N\in\naturals$ such that every
  $I\in H_f(R)$ has $\dim_\FF (R/I)_N = f(N)$, and such that the map
  $I \mapsto (R/I)_N \in Gr^{f(N)}(R_N)$ is a closed embedding of
  $H_f(R)$ into the Grassmannian of $f(N)$-dimensional quotients of $R_N$.
\end{Theorem*}

The tangent space at $I$ on this Grassmannian is $Hom_\FF(I_N,(R/I)_N)$.
The tangent space to the Hilbert scheme is defined by the $R$-linear,
not just $\FF$-linear, maps \cite[Thm. IV-19]{WhySchemes}.

The subscheme $H_f(R,\varphi)$ will be the intersection of $H_f(R)$
with another subscheme $\Phi \subseteq Gr^{f(N)}(R_N)$. 
Let
$$ 
\Phi = \left\{ V \in Gr^{f(N)}(R_N): \forall r\in V, \varphi(r)^p \in V \right\}
$$
(where $R_N$ is thought of as an $\FF$-vector space with respect to
the twisted action $c\cdot r := c^p r$, to make the these equations algebraic),
and let $H_f(R,\varphi) = H_f(R) \cap \Phi$. 

\begin{Theorem}\label{thm:main}
  Let $R$ be a graded ring with $R_0 = \FF$ and a graded splitting.
  Let $I \in H_f(R,\varphi)$ be a graded, compatibly split ideal. 
  Then the Zariski tangent space to $H_f(R,\varphi)$ at $I$ vanishes.  
  In particular $H_f(R,\varphi)$ is a finite set of reduced points.
\end{Theorem}

\begin{proof}
  As in \cite[Thm. IV-19]{WhySchemes},
  it is easy to see that the tangent space consists of $R$-linear maps
  $\lambda: I\to R/I$ intertwining the restriction $\varphi_I$ of the
  splitting $\varphi$ on $R$ and the induced splitting $\varphi_{R/I}$ on $R/I$.

  Now let $i\in I$ and let $\lambda : I \to R/I$
  intertwine $\varphi_I,\varphi_{R/I}$. Then
  \begin{align*}
       \lambda(i)
    &= \lambda(\varphi_I(i^p))&\text{since $\varphi_I$ is restricted from $R$}\\
    &= \varphi_{R/I}(\lambda(i^p))  &\text{the intertwining condition}\\
    &= \varphi_{R/I}(i^{p-1} \lambda(i)) &\text{since $\lambda$ is $R$-linear}\\
    &= \varphi_{R/I}(0)             &\text{since on $R/I$, $i$ acts as $0$,
    and $p-1 \geq 1$}\\
    &= 0 
  \end{align*}
  so for any $\lambda$ in the tangent space, $\lambda = 0$.

  Since $H_f(R,\varphi)$ is of finite type, it has finitely many components,
  each of which is a reduced point by the above calculation.
\end{proof}

This does not quite recover Schwede's theorem (even in the graded case);
while it shows that there are only finitely many compatibly split 
subschemes with given Hilbert polynomial, it doesn't show that only
finitely many Hilbert polynomials arise.

Either theorem is enough to prove the following

\begin{Corollary}\label{cor:invt}
  Let $R$'s grading extend to a multigrading, and assume the splitting $\varphi$
  is multigraded (perhaps by replacing $\varphi$ by its multigraded part).
  Then any compatibly split saturated singly graded ideal $I$ is
  automatically multigraded.

  Using Schwede's theorem, one need not even assume from the outset
  that $I$ is singly graded.
\end{Corollary}

\begin{proof}
  Rather than the multigrading, we work with the equivalent action 
  of a torus $T$ on $R$. Then since $\varphi$ is $T$-invariant, 
  $T$ acts on the set of compatibly split ideals.
  Since this set is discrete and $T$ is connected, it fixes each ideal.
\end{proof}

\section{Examples}

\newcommand\Proj{{\rm Proj}}

In this section we use geometric language, with $X = \Proj\ R$,
and speak of subschemes rather than saturated graded ideals.
We continue to interpret multigradings as torus actions,
and multigraded splittings as torus-invariant.

\subsection{Compatibly split points}

Any single point in $X$ has the Hilbert polynomial $f\equiv 1$, so
the set of compatibly split points is the finite set $H_{1}(R,\varphi)$.

If $X$ carries a torus action and $\varphi$ is $T$-invariant, then 
Corollary \ref{cor:invt} tells us
$H_{1}(R) \subseteq X^T$, where $X^T$ is the fixed-point set.

I do not know if isolated $T$-fixed points are automatically compatibly split.
It is not hard to prove that if $F$ is an attractive component of $X^T$
(meaning, the sink in some Bia\l ynicki-Birula decomposition of $X$)
then $F$ is compatibly split. So for example, any $T$-invariant
splitting on a flag manifold (as classified in \cite[Chapter 2.3]{BrionKumar})
splits all the $T$-fixed points, since they are all attractive.

\subsection{$X$ a toric variety}

If $X$ carries a torus action with a dense orbit, then its normalization 
is a toric variety, and it has only finitely many orbits. Hence $X$ has only
finitely many reduced $T$-invariant subschemes.
If $X$ has a splitting, then it has a $T$-invariant splitting, 
which turns out to be unique. With respect to this splitting,
any compatibly split subscheme is reduced and $T$-invariant,
so we recover Schwede's theorem in this simple case.

Being split does not imply that $X$ is normal; rather, it is equivalent
(in the presence of the dense torus orbit assumption) to $X$ being seminormal.
For example, let $X = {\mathbb P}^2 / \sim$
where $\left\{ [a,b,0] \sim [-a,b,0] \right\}$ generate the equivalence 
relation $\sim$; this $X$ is Frobenius split and seminormal but 
not normal (its normalization is ${\mathbb P}^2$).

\subsection*{Acknowledgements} 
We thank Michel Brion, Shrawan Kumar, and Karl Schwede for
useful discussions.

\end{document}